\documentclass{amsart}
\usepackage{amssymb}
\usepackage{url}

\newtheorem{lemma}{Lemma}

\newcommand{\R}{I\!R}             
\newcommand{\N}{I\!N}             
\newcommand{\deviation}{\ensuremath{c}}
\newcommand{\demand}{\ensuremath{d}}
\newcommand{\maxlots}{\ensuremath{k}}

\begin{document}

\title{The combinatorics of S, M, L, XL:\\The best fitting delivery of T-shirts}

\author{Constantin Gaul, Sascha Kurz$^{\ast}$ and J\"org Rambau}
\address{firstname.lastname@uni-bayreuth.de, University of Bayreuth, Dep.~of Mathematics, Physics, and Computer Science\\
  95440 Bayreuth, Germany}

\maketitle

\begin{abstract}
  We consider the problem of approximating the branch and size
  dependent demand of a fashion discounter with many branches by a
  distributing process being based on the branch delivery restricted
  to integral multiples of lots from a small set of available
  lot-types. We propose a formalized model which arises from a
  practical cooperation with an industry partner. Besides an integer
  linear programming formulation we provide an appropriate primal
  heuristic for this problem.

  \bigskip

\end{abstract}

\thanks{This research was supported by a grant of the \emph{Bayerische
  Forschungsstiftung, Az-727-06}.}

\section{Introduction \label{intro}}
The problem studied in this note is motivated by a special feature
of the ordering process of a fashion discounter with many branches:
For each product that hits the shelves, the internal stock-turnover
has to distribute around 10\,000 pieces among the around 1\,000
branches, correctly assorted by size.  This would mean 10\,000 picks
with high error probability in the central-warehouse (in our case in
the high-wage country Germany). In order to reduce the handling costs
and the error proneness in the central warehouse, all products are
ordered in multiples of so-called \emph{lot-types} from the suppliers
who in general are located in extremely low-wage countries.

A lot-type specifies a number of pieces of a product for each
available size, e.g., (2,2,2,2,2) if the sizes are (S, M, L, XL, XXL)
means two pieces of each size.  A \emph{lot} of a certain lot-type is
a foiled pre-pack that contains as many pieces of each size as
specified in its lot-type. The number of different lot-types is
bounded by the supplier.

So we face an approximation problem: which (integral) multiples of
which (integral) lot-types should be supplied to a branch in order to
meet a (fractional) mean demand as closely as possible?  We call this
specific demand approximation problem the \emph{lot-type design
  problem (LDP)}.  A detailed version of this work appeared
in~\cite{Gaul+Kurz+Rambau:LotTypeDesignProblemOMS:2009}, where also
references to related work can be found.


\section{The lot-type design problem\label{sec_lottype}}

Formally, the problem can be stated as follows: Consider a fashion
discounter with branches $b \in \mathcal{B}$ who wants to place an
order for a certain product that can be obtained in sizes $s \in
\mathcal{S}$ and that can be pre-packed in lot-types $l \in
\mathcal{L}$.  Each lot-type is a vector $(l_s)_{s \in \mathcal{S}}$
specifying the number of pieces of each size contained in the
pre-pack.  Only $\maxlots$ different lot-types from $\mathcal{L}$ are
allowed in this order, and each branch receives only lots of a single
lot-type.  We are given lower and upper bounds $\underline{I},
\overline{I}$ for the total supply of this product.  Moreover, we
assume that a the branch and size dependent mean demand $\demand_{b,
  s}$ for the corresponding type of product is known to us.

The original goal is to find a set of at most $\maxlots$ lot-types, an
order volume for each of these chosen lot-types, and a distribution of
lots to branches such that the revenue is maximized.  In order to
separate the order process from the sales process (which involves
mark-downs, promotions, etc.), we restrict ourselves in this paper to
the minimization of the distance between supply and mean demand
defined by a vector norm.

The \emph{Lot-Type Design Problem (LDP)} is the following
optimization problem:

\begin{center}
  \begin{tabular}[t]{rp{0.8\linewidth}}
    \emph{Instance:} & 
    We are given
    \begin{itemize}
    \item a set of branches $b \in \mathcal{B}$
    \item a set of sizes $s \in \mathcal{S}$
    \item a mean demand table $\demand_{b, s}$, $b \in \mathcal{B}$, $s
      \in \mathcal{S}$
    \item a norm $\lVert\cdot\rVert$ on~$\R^{\mathcal{B}
        \times \mathcal{S}}$
    \item a set $\mathcal{L}$ of feasible lot types $(l_s)_{s \in
        \mathcal{S}} \in \N_0^{\mathcal{S}}$
    \item a maximal number $M \in \N$ of possible
      multiplicities
    \item a maximal number $\maxlots \in \N$ of lot types to
      use
    \item lower and upper bounds $\underline{I}$, $\overline{I}$ for
      the total supply
    \end{itemize}\\
    \emph{Task:} & For each branch~$b \in \mathcal{B}$ choose a lot
    type $l(b) \in \mathcal{L}$ and a number $m(b) \in \N$,
    $1 \le m(b) \le M$ of
    lots to order for~$b$ such that
    \begin{itemize}
    \item the total number of ordered lot types is at most~$\maxlots$
    \item the total number of ordered pieces is in $[\underline{I},
      \overline{I}]$\newline(the \emph{total capacity constraint})
    \item the distance of the order from the demand measured by
      $\lVert\cdot\rVert$ is minimal
    \end{itemize}
  \end{tabular}
\end{center}

The LDP can be formulated as an Integer Linear Program if we restrict
ourselves to the $L^1$-norm for measuring the distance between supply
and demand.  This norm is quite robust against outlies in the demand
estimation.

We use binary variables $x_{b,l,m}$, which are equal to $1$ if and
only if lot-type $l$ is delivered with multiplicity $m$ to Branch $b$,
and binary variables $y_l$, which are $1$ if and only if at least one
branch in $\mathcal{B}$ is supplied with Lottype~$l$.  The
\emph{deviation} of the demand from the supply if Branch~$b$ is
supplied by $m$ lots of lot-type~$l$ is given by $\deviation_{b,l,m}
:= \sum_{s \in \mathcal{S}} \lvert \demand_{b, s} - m \cdot l_s \rvert$.

The following integer linear program models the LDP with $L^1$-norm.
\begin{align}
  \label{OrderModel_Target}
  \min && \sum_{b\in\mathcal{B}}\sum_{l\in\mathcal{L}}\sum_{m=1}^M \deviation_{b,l,m}\cdot x_{b,l,m}\\
  \label{OrderModel_EveryBranchOneLottype}
  s.t. 
  && 
  \sum_{l\in\mathcal{L}}\sum_{m=1}^M x_{b,l,m} &= 1 && \forall b\in\mathcal{B}\\
  \label{OrderModel_UsedLottypes}
  && 
  \sum_{l\in\mathcal{L}} y_l & \le \maxlots\\
  \label{OrderModel_Binding}
  && 
  \sum_{m=1}^M x_{b,l,m} & \le y_l && \forall b\in\mathcal{B}, \forall
  l\in\mathcal{L}\\
  &&
  \label{OrderModel_Cardinality}
  \underline{I} \le \sum_{b\in\mathcal{B}}\sum_{l\in\mathcal{L}}\sum_{m=1}^M
  \sum_{s \in \mathcal{S}} 
  m \cdot l_s \cdot x_{b,l,m} &\le \overline{I}\\
  && 
  x_{b,l,m} & \in\{0,1\} && \forall b\in\mathcal{B}, \forall
  l\in\mathcal{L}, \forall m = 1,\dots,M\\
  && 
  y_l & \in\{0,1\} && \forall l\in\mathcal{L}
\end{align}

The objective function \eqref{OrderModel_Target} computes the $L^1$-distance
of the supply specified by $x$ from the demand.  Condition
\eqref{OrderModel_EveryBranchOneLottype} enforces that each branch is
assigned a unique lot-type and a unique multiplicity.  Condition
\eqref{OrderModel_UsedLottypes} models that at most $\maxlots$ different
lot-types can be chosen.  Condition \eqref{OrderModel_Binding} forces
the selection of a lot-type whenever it is assigned to some branch
with some multiplicity.  Finally, Condition
\eqref{OrderModel_Cardinality} ensures that the total number of pieces
is in the desired interval -- the total capacity constraint.

Our ILP formulation can be used to solve all real world instances of
our business partner in at most 30~minutes by using a standard ILP
solver like \texttt{CPLEX 11}.  Interestingly, the model seems quite
tight -- most of the time is spent in solving the root LP.  

Although 30~minutes may mean a feasible computation time for an
offline-optimization in many contexts, this is not fast enough for our
real world application. The buyers of our retailer need a software
tool which can produce a near optimal order recommendation in real
time on a standard laptop.  For this reason, we present a fast
anytime-heuristic, which has only a small gap compared to the optimal
solution on a test set of real world data of our business partner.

We briefly sketch the idea of the heuristic \emph{Score-Fix-Adjust
  (SFA)}: It
\begin{enumerate}
\item sorts all lot-types according to certain scores, coming from a
  count for how many branches the lot-type fits best, second best,
  \ldots (Score);
\item fixes $\maxlots$-subsets of lot-types in the order of decreasing
  score sums (Fix);
\item greedily adjusts the multiplicities so as to achieve feasibility
  w.r.t.\ the total capacity
  constraint (Adjust).
\end{enumerate}
Details can be found in \cite{Gaul+Kurz+Rambau:LotTypeDesignProblemOMS:2009}.

Since in the case $\maxlots = 1$ we can very often loop over all
feasible lot-types, it is interesting that in this case SFA always
yields an optimal solution (for any norm).

\begin{lemma}
  For $\maxlots = 1$ and costs $c_{b,l,m}=\Vert d_{b,\cdot}-m\cdot l\Vert$
  for an arbitrary norm
  $\left\Vert\cdot\right\Vert$, our heuristic SFA produces an optimal
  solution whenever all lot-types $l\in\mathcal{L}$ are checked.
\end{lemma}

In order to substantiate the usefullness of our heuristic, we have
compared the quality of the solutions, given by this heuristic after
one second of computation time (on a standard laptop:
Intel$^{\textregistered}$ Core$\texttrademark\ $2 CPU with 2~GHz and
1~GB RAM) with respect to the solution given by \texttt{CPLEX 11}
(after solving to optimality).

Our business partner has provided us with historic sales information
for nine different commodity groups, each ranging over a sales period
of at least one-and-a-half years.  From this we estimated mean demands
via aggregating over products in a commodity group.  By normalizing
the lengths of the products' sales periods to the point in time when
half of the product was sold out, we were able to mod out the effects
of any product's individual success or failure.  Prior to each test
calculation, the resulting demands were scaled so that the total mean
demand was in the center of the total capacity interval given by the
management for a new order of a product in that commodity group.

For each commodity group we have performed a test calculation for
$\maxlots\in\{2,3,4,5\}$ distributing some amount of items to almost
all branches. The crucial parameters are given in Table
\ref{table_parameters}, the results are presented in
Table~\ref{table_gap}.

\begin{table}[htp]\footnotesize\sffamily
  \begin{center}\renewcommand{\arraystretch}{1.41}
    \begin{tabular}{r@{\hspace*{1cm}}c@{\hspace*{1cm}}c@{\hspace*{1cm}}c@{\hspace*{1cm}}c@{\hspace*{1cm}}c}
      \hline
      Commodity group & $|\mathcal{B}|$ & $|\mathcal{S}|$ & $\left[\underline{I},\overline{I}\right]$
      & $|\mathcal{L}|$ & $M$\\
      \hline 
      1 & 1119 & 5 & [10\,630, 11\,749] & 243 & 10 \\
      2 & 1091 & 5 & [10\,000, 12\,000] & 243 & 10 \\
      3 & 1030 & 5 & [9\,785, 10\,815]  & 243 & 10 \\
      4 & 1119 & 5 & [10\,573, 11\,686] & 243 &  9 \\
      5 & 1175 & 5 & [16\,744, 18\,506] & 243 & 15 \\
      6 & 1030 & 5 & [11\,000, 13\,000] & 243 &  9 \\
      7 & 1098 & 5 & [15\,646, 17\,293] & 243 &  9 \\
      8 &  989 & 5 & [11\,274, 12\,461] & 243 &  9 \\
      9 &  808 & 5 & [9\,211, 10\,181]  & 243 & 10 \\
      \hline
    \end{tabular}
    \caption{Parameters for the test calculations.}
    \label{table_parameters}
  \end{center}
\end{table}

\begin{table}[htp]\footnotesize\sffamily
  \begin{center}\renewcommand{\arraystretch}{1.5}
    \begin{tabular}{r@{\hspace*{1cm}}c@{\hspace*{1cm}}c@{\hspace*{1cm}}c@{\hspace*{1cm}}c}
      \hline
      Commodity group & $\maxlots=2$ & $\maxlots=3$ & $\maxlots=4$ & $\maxlots=5$ \\
      \hline
      1 & 2.114\,\% & 1.226\,\% & 2.028\,\% & 1.983\,\% \\
      2 & 0.063\,\% & 0.052\,\% & 0.006\,\% & 0.741\,\% \\
      3 & 0.054\,\% & 0.094\,\% & 0.160\,\% & 0.170\,\% \\
      4 & 0.019\,\% & 0.007\,\% & 0.024\,\% & 0.038\,\% \\
      5 & 0.015\,\% & 0.017\,\% & 0.018\,\% & 0.019\,\% \\
      6 & 0.018\,\% & 0.022\,\% & 0.024\,\% & 0.022\,\% \\
      7 & 0.013\,\% & 0.013\,\% & 0.014\,\% & 0.014\,\% \\
      8 & 0.016\,\% & 0.017\,\% & 0.018\,\% & 0.019\,\% \\
      9 & 0.011\,\% & 0.939\,\% & 0.817\,\% & 0.955\,\% \\
      \hline
    \end{tabular}
    \caption{Optimality gap in the $\Vert\cdot\Vert_1$-norm for our heuristic on nine commodity groups and
      different values for the maximum number $\maxlots$ of used lot-types.}
    \label{table_gap}
  \end{center}
\end{table}

We can see that -- given the uncertainty in the data -- the
performance of SFA is more than satisfactory.

\section{Conclusions}

We identified the lot-type design problem in the supply chain
management of a fashion discounter.  It can be modeled as an ILP, and
real-world instances can be solved by commercial-of-the-shelf software
like CPLEX in half an hour whenever the number of lot-types is not too large.

Our SFA-heuristics finds solutions with a gap of mostly under 1\,\% in
a second, also for instances with a large number of lot-types.  Given
the volatility of the demand estimation, these gaps are certainly
tolerable.

Meanwhile, the model and SFA have been put to operation by our
business partner with significant positive monetary impact.


\end{document}